\documentclass[12pt]{amsart}

\usepackage{latexsym}
\usepackage{amsmath}
\usepackage{amssymb}
\usepackage{amscd}
\usepackage{multicol}

\newtheorem{theorem}{\bf Theorem}[section]
\newtheorem{corollary}[theorem]{\bf Corollary}
\newtheorem{proposition}[theorem]{\bf Proposition}
\newtheorem{lemma}[theorem]{\bf Lemma}

\newtheorem{definition-theorem}[theorem]{\bf Theorem-Definition}

\def\h{{\rm ht}}

\def\bR{\mathbb{R}}
\def\bC{\mathbb{C}}

\def\t{\mathfrak{t}}

\def\la{(\! (}
\def\ra{)\! )}

\def\t{\frak{t}}

\def\H{{\mathcal H}}
\title[Combinatorial quantum cohomology of $G/B$]{The combinatorial
quantum cohomology ring of $G/B$}

\begin{document}
\begin{abstract}
A purely combinatorial construction of the quantum cohomology ring
of the generalized flag manifold is presented. We show that the
ring we construct is commutative, associative and satisfies the
usual grading condition. By using results of our previous papers
[Ma1] and [Ma2], we obtain a presentation of this ring in terms of
generators and relations, and formulas for quantum Giambelli
polynomials. We show that these polynomials satisfy a certain
orthogonality property, which --- for $G=SL_n(\bC)$ --- was proved
previously in the paper [Fo-Ge-Po].
\end{abstract}

\author[A.-L. Mare]{Augustin-Liviu  Mare}
\address{Department of Mathematics\\ University of Toronto
 \\Toronto, Ontario M5S 3G3, Canada}
 \email{amare@math.toronto.edu}

\maketitle

\section{Introduction}

Let us consider the complex flag manifold $G/B$, where $G$ is a
connected, simply connected,  simple, complex Lie group and
$B\subset G$ a Borel subgroup. Let  $\t$ be the Lie algebra of a
maximal torus of a compact real form of $G$ and  $\Phi \subset
\t^*$ the corresponding set of roots. Consider an arbitrary
 $W$-invariant inner product $\langle \ , \  \rangle$ on $\t$.
 The Weyl group $W$ is
the subgroup of $O(\t, \langle \ , \  \rangle)$ generated by the
reflections about the hyperplanes $\ker \alpha $,
 $\alpha\in \Phi^+$.
 To any root $\alpha$ corresponds the coroot
$$\alpha^{\vee}:= \frac{2\alpha}{\langle \alpha,
\alpha\rangle}$$ which  is an element of $\t$, by using the
identification of $\t$ and $\t^*$ induced by $\langle \ , \
\rangle$. If $\{\alpha_1, \ldots ,\alpha_l\}$ is a system of
simple roots then $\{\alpha_1^{\vee},\ldots, \alpha_l^{\vee}\}$ is
a  system of simple coroots.  Consider  $\{\lambda_1 ,\ldots ,
\lambda_l\} \subset \t^*$ the corresponding system of fundamental
weights, which are defined by
$\lambda_i(\alpha_j^{\vee})=\delta_{ij}$. It can be shown that the Weyl
group $W$ is actually generated by the {\it simple reflections} $s_1=s_{\alpha_1},\ldots,
s_l=s_{\alpha_l}$ about the hyperplanes $\ker\alpha_1, \ldots, \ker\alpha_l$.
The {\it length} $l(w)$ of $w$ is the minimal number of factors in a decomposition
of $w$ as a product of simple reflections. We denote by $w_0$ the longest element
of $W$.

Let $B^-\subset G$ denote the Borel subgroup opposite to $B$. To
each $w\in W$ we assign the {\it Schubert variety}
$X_w=\overline{B^-.w}$. The Poincar\'e dual of $[X_w]$ is an
element of $H^{2l(w)}(G/B)$, which is called the {\it Schubert
class}. The set $\{\sigma_w~|~w\in W\}$ is a basis of the
cohomology\footnote{Only cohomology with coefficients in $\bR$
will be considered in this paper.} module $H^*(G/B)$. The
Poincar\'e pairing $( \ , \ )$ on $H^*(G/B)$ is determined by:
\begin{equation}\label{poincare} ( \sigma_u,\sigma_v )
= \left\{
\begin{array}{ll}
                         1,\ {\rm if }\ u=w_0v  \\
                         0, \ {\rm otherwise}
                             \end{array}
                      \right.\end{equation}

According to  a theorem of Borel [Bo], the  ring homomorphism
$S(\t^*)\to H^*(G/B)$ defined by $$\lambda_i \mapsto \sigma_{s_i},
\quad 1\le i \le l$$ is surjective and it induces the  ring
isomorphism
\begin{equation}\label{borel} H^*(G/B) \simeq \bR[\{\lambda_i\}]/I_W, \end{equation}
where $I_W$ is the ideal of
$S(\t^*)=\bR[\lambda_1,\ldots,\lambda_l]=\bR[\{\lambda_i\}]$
generated by the $W$-invariant polynomials of strictly positive
degree. Recall that, by a result of Chevalley [Ch],  there exist
$l$ homogeneous, functionally independent  polynomials
$u_1,\ldots, u_l \in S(\t^*)$ which generate $I_W$.  We identify
$H^*(G/B)$ with Borel's presentation and denote them both by $
\H$. So
$$\H=H^*(G/B)=\bR[\{\lambda_i\}]/I_W.$$ In this way the homogeneous elements of $\H$ will be
of the form$$[f]=f \ {\rm  mod} \  I_W,$$ where $f\in
\bR[\{\lambda_i\}]$ is a homogeneous polynomial. In particular,
the degree two Schubert classes will be $[\lambda_i]$, $1\le i \le
l$.

In fact we would like to see {\it all}
 Schubert classes as cosets of
certain polynomials  in the presentation  (\ref{borel}). A
construction of such polynomials  was obtained by Bernstein, I. M.
Gelfand and S. I. Gelfand in [Be-Ge-Ge], as follows: To each
positive root $\alpha$ we assign the {\it divided difference
operator} $\Delta_{\alpha}$ on the ring $\bR[\{\lambda_i\}]$
(since the latter is just the symmetric algebra $S(\t^*)$,  it
admits a natural action of the Weyl group $W$):
$$\Delta_{\alpha}(f)=\frac{f-s_{\alpha}f}{\alpha},$$
$f\in \bR[\{\lambda_i\}]$.  If $w$ is an arbitrary element of
$W$, take $w =s_{i_1} \ldots s_{i_k}$  a reduced expression and
then set\footnote{ One can show (see for instance [Hi, Ch. IV])
that the definition does not depend on the choice of the reduced
expression.} $$\Delta _w= \Delta _{\alpha _{i_1}} \circ \cdots
\circ \Delta_{\alpha _{i_k}}.$$ The polynomial
$$c_{w_0}:=\frac{1}{|W|}\prod_{\alpha \in \Phi^+}\alpha$$ is
homogeneous, of degree $l(w_0)$ and
 has the  property that $\Delta_{w_0}c_{w_0}=1$.
 To any $w\in W$ we assign $$c_w:= \Delta_{w^{-1}w_0}c_{w_0}$$
which is a homogeneous polynomial  of degree $l(w)$ satisfying
\begin{equation}\label{del}\Delta_v c_w= \left\{
\begin{array}{ll}
   c_{wv^{-1}}, & {\rm if \ } l(wv^{-1})=l(w)-l(v) \\
    0, & {\rm otherwise}
\end{array}\right.
\end{equation} for any $v\in W$ (see for instance [Hi, Ch. IV]).

 \begin{theorem}{\rm ([Be-Ge-Ge])}\label{bernstein} By the identification (\ref{borel}) we have
$$\sigma_w=[c_w],$$
for any $w\in W$.
\end{theorem}

The main goal of our paper is to construct in a purely
combinatorial way a certain ``quantum deformation'' of the ring
$\H$. This will depend on the ``deformation parameters''
$q_1,\ldots,q_l$, which are just some additional multiplicative
variables. Let us begin with the following lemma, which was proved
for instance in [Ma1] (see also \cite{Pe} or \cite{Br-Fo-Po}).
Recall first that if $\alpha$ is a positive root, then the {\it
height} of the corresponding coroot $\alpha^{\vee}$ is by
definition
$$\h (\alpha^{\vee})=m_1+\ldots +m_l,$$ where the positive integers
$m_1,\ldots, m_l$ are given by
\begin{equation}\label{length}\alpha^{\vee}=m_1\alpha_1^{\vee}+\ldots +m_l\alpha_l^{\vee}.
\end{equation}
\begin{lemma}\label{inequality}
 For any positive root $\alpha$ we have that
$l(s_{\alpha})\leq 2 \h(\alpha^{\vee})-1$.
\end{lemma}
Denote by $\tilde{\Phi}^+$ the set  of all positive  roots
$\alpha$ with the property that
$$l(s_{\alpha})=2\h (\alpha^{\vee})-1.$$
We will obtain in section 3 a complete description of the elements
of $\tilde{\Phi}^+$ (see Lemma 3.1). One can easily deduce from
this that if the root system of $G$ is simply laced, then
$\tilde{\Phi}^+ =\Phi^+$.

 The following divided
difference type operators on $\bR[\{\lambda_i\},\{q_i\}]$ have
been considered by D. Peterson in [Pe]:
\begin{equation}\label{Lambda}\Lambda_j=\lambda_j +\sum_{\alpha\in
\tilde{\Phi}^+}\lambda_j
(\alpha^{\vee})q^{\alpha^{\vee}}\Delta_{s_{\alpha}}, \quad 1\le j
\le l\end{equation} where we use the notation
$$q^{\alpha^{\vee}}=q_1^{m_1}\ldots q_l^{m_l},$$
with $m_1\ldots,m_l$ given by (\ref{length}). It is obvious that
$\Lambda_j$ leaves the ideal $I_W\otimes \bR[\{q_i\}]$ of
$\bR[\{\lambda_i\},\{q_i\}]$ invariant, hence it induces an
operator on $\H\otimes \bR[\{q_i\}]$.

The following result\footnote{The proof of this result given in
[Ma1] relies essentially on the associativity of the ring
$QH^*(G/B)$, which  is a highly nontrivial fact; the proof of
Lemma \ref{peterson} we are going to give here is entirely in the
realm of root systems.} was stated by D. Peterson in [Pe] (for
$G=SL(n,\bC)$, a proof can be found in \cite{Fo-Ge-Po}).
\begin{lemma}\label{peterson} The operators $\Lambda_1,\ldots,\Lambda_l$ on
$\bR[\{\lambda_i\},\{q_i\}]$ commute.
\end{lemma}
We will prove this lemma in section 3 of our paper.  The operator
$\psi$ defined in the next lemma will be an important object in
our paper.
\begin{lemma}\label{thirdlemma} The map $\psi:\bR[\{\lambda_i\},\{q_i\}]\to \bR[\{\lambda_i\},\{q_i\}]$
given by
$$\psi(f)=f(\{\Lambda_i\},\{q_i\})(1),$$
$f\in \bR[\{\lambda_i\},\{q_i\}]$ is an isomorphism of
$\bR[\{q_i\}]$-modules. For $f\in \bR[\{\lambda_i\},\{q_i\}]$ of
degree $m$ with respect to $\lambda_1,\ldots ,\lambda_l$, we have
\begin{align}
\psi^{-1}(f)=&\frac{I-(I-\psi)^{m}}{\psi}(f)\nonumber
\\=&{m\choose 1} f-{{m}\choose 2}\psi(f) +\ldots + (-1)^{m
-2}{{m}\choose{m-1}}\psi^{m -2}(f)+(-1)^{m-1}\psi^{m
-1}(f),\nonumber  \end{align} where ${m\choose 1} ,\ldots , {m
\choose {m-1}}$ are the binomial coefficients.
\end{lemma}
The proof follows in an elementary way from  the fact that the
degree of $f-\psi(f)$ with respect to $\lambda_1,\ldots
,\lambda_l$  is strictly less than the degree of $f$ (the details
can be found in [Ma1, Lemma 3.4]).

Our aim is to investigate the ring defined as follows (note that
for $G=SL(n,\bC)$ a similar object has been considered by A.
Postnikov in [Po]).
\begin{definition-theorem}\label{defthm} The composition law $\star$ on
the $\bR[\{q_i\}]$-module $\H\otimes
\bR[\{q_i\}]=\bR[\{\lambda_i\},\{q_i\}]/(I_W\otimes \bR[\{q_i\}])$
given by
\begin{equation}\label{defstar}[f]\star
[g]=[\psi(\psi^{-1}(f)\psi^{-1}(g))],\quad f,g \in
\bR[\{\lambda_i\},\{q_i\}]\end{equation} is well defined,
commutative, associative, $\bR[\{q_i\}]$-bilinear, and  satisfies:
\begin{itemize}
\item
$\deg (a\star b)=\deg a +\deg b$, for any two homogeneous elements
$a, b$ of $\H\otimes \bR[\{q_i\}]$, where we assign $$\deg
[\lambda_i]=2, \ \deg q_i=4,\quad 1\le i \le l.$$
\item (Frobenius property) $( a\star b,c) =(
a, b\star c)$, for any $a,b,c\in \H$, where $(  \ , \  )$ is the
$\bR[\{q_i\}]$-bilinear extension of the Poincar\'e pairing on
$\H$.
\end{itemize}

 We will call
$\star$ the {\rm combinatorial quantum  product} on
$\H\otimes \bR[\{q_i\}]$.
\end{definition-theorem}

We will prove this theorem at the beginning of section 2.

A complete knowledge of the combinatorial quantum cohomology
$\bR[\{q_i\}]$-algebra
 defined in the previous theorem can be achieved by finding the
structure constants (which are in $\bR[\{q_i\}]$) of the
multiplication $\star$ with respect to the basis consisting of the
Schubert classes $\sigma_w = [c_w]$, $w\in W$. Like in the
classical situation (see the beginning of this section), we can
obtain this information about $(\H\otimes \bR[\{q_i\}],\star)$ as
follows:
\begin{itemize}
\item[(a)] describe it in terms of generators and relations (i.e.
find the quantum analogue of Borel's presentation (\ref{borel}))
\item[(b)] determine representatives of the Schubert classes in
the quotient ring obtained at (a) (i.e. find the quantum analogue
of the Bernstein-Gelfand-Gelfand polynomials, see Theorem
\ref{bernstein}).
\end{itemize}

The next two theorems give solutions to problems (a), respectively
(b). The first theorem can be interpreted as the combinatorial
version of B. Kim's theorem [Kim]. Our proof, which can be found
in section 2, is a direct application of a more general result
obtained by us in [Ma2].

\begin{theorem}\label{firsttheorem} Let   $I_W^q$ denote the ideal of
 $\bR[\{ \lambda_i\},\{q_i\}]$ generated by
 $F_k(\{\lambda_i\},\{ -\langle \alpha_i^{\vee},\alpha_i^{\vee}\rangle q_i\} )$,
$1\leq k \leq l$, where  $F_k$  are polynomials in $2l$ variables
which represent the integrals of motion of the Hamiltonian system
of Toda lattice type associated to the coroot system of $G$ (for
more details, see section 2). Then the map
$$(\H \otimes
\bR[\{q_i\}] = \bR[\{\lambda_i\},\{q_i\}]/(I_W\otimes \bR[q_i]),
\star) \to\bR[\{ \lambda_i\},\{q_i\}]/I_W^q, $$ given by $$f \
{\rm mod} \ I_W  \ \mapsto \psi^{-1}(f) \ {\rm mod} \ I_W^q,$$
$f\in \bR[\{\lambda_i\},\{q_i\}]$, is an isomorphism of
$\bR[\{q_i\}]$-algebras.
\end{theorem}

Alternatively, one can see that $I_W^q$ is the ideal of
$\bR[\{\lambda_i\},\{q_i\}]$ generated by the polynomials
$\psi^{-1}(u_1),\ldots, \psi^{-1}(u_l)$, which is the same as
$\psi^{-1}(I_W)$ (see Proposition \ref{psiminus}).

What follows now is the combinatorial version of the main result
of [Ma1], where a quantum Giambelli formula for $G/B$ has been
obtained. In the context of our present paper, we obtain the same
formula by a straightforward application of Theorem
\ref{firsttheorem} and Lemma \ref{thirdlemma}.

\begin{corollary}\label{secondtheorem} The  isomorphism described
by Theorem \ref{firsttheorem} maps the Schubert class
$\sigma_w=c_w \  {\rm mod} \ I_W $ to the class modulo $I_W^q$ of
the polynomial
\begin{align}
\psi^{-1}(c_w)=&\frac{I-(I-\psi)^{l}}{\psi}(c_w)\nonumber
\\=&{l\choose 1} c_w-{{l}\choose
2}\psi(c_w) +\ldots + (-1)^{l -2}{{l}\choose{l-1}}\psi^{l
-2}(c_w)+(-1)^{l-1}\psi^{l -1}(c_w),\nonumber  \end{align}
 where $l$ denotes $l(w)$.
\end{corollary}

We will also show that the polynomials described by Corollary
\ref{secondtheorem} satisfy a certain orthogonality condition
(similar to (\ref{poincare})) with respect to the ``quantum
intersection pairing'' (see Proposition \ref{prop}) .

{\bf Remarks.} 1. The actual  quantum product $\circ$ on
$H^*(G/B)\otimes \bR[\{q_i\}]$  is defined in terms of numbers of
holomorphic curves which intersect ``general'' translates of three
given Schubert varieties (for the precise definition, one can see
\cite{Fu-Pa} or \cite{Fu-Wo}). The {\it quantum Chevalley formula}
describes the  multiplication of degree two Schubert classes by
arbitrary Schubert classes. More precisely, in terms of the
identification (\ref{borel}) (see also Theorem \ref{bernstein}),
it states that
\begin{equation}\label{quantumc}[\lambda_i]\circ [c_w]=\Lambda_i ([c_w]).\end{equation}
This formula was announced by Peterson in \cite{Pe} and then
proved by Fulton and Woodward in \cite{Fu-Wo}. In order to relate
(\ref{quantumc}) to our product $\star$, we  note that
\begin{equation}\label{Lambdastar} \Lambda_i
([c_w])=[\Lambda_i(c_w)]= [\Lambda_i(\psi\psi^{-1}({c}_w))]=
[\psi(\lambda_i\psi^{-1}({c}_w))]=[\lambda_i]\star[c_w]\end{equation}
where we have used that $\psi(\lambda_i)=\lambda_i$. We deduce
that
$$[\lambda_i]\circ[c_w]=[\lambda_i]\star[c_w], \quad 1\le i \le l,
w\in W.$$ This implies that
$$[c_v]\circ [c_w]=[c_v]\star[c_w],$$ for any $v,w\in W$, because
both $(\H\otimes \bR[\{q_i\}], \star)$ and $(\H\otimes
\bR[\{q_i\}], \circ)$ are generated by $[\lambda_1],\ldots,
[\lambda_l]$ as $\bR[\{q_i\}]$-algebras.  Now, since
$\star=\circ$, the results about $\star$ which we  prove in our
paper
 hold for $\circ$ as well. In this way we are able to prove results about the
 actual quantum cohomology ring
$QH^*(G/B)=(\H\otimes \bR[\{q_i\}],\circ).$ Except for Proposition
\ref{prop}
--- which was proved in [Fo-Ge-Po] for $G=SL(n,\bC)$  --- these
results are not new (see  [Kim]   and [Ma1]).

2. We hope that that a similar  approach can be used by
considering instead of the root system of $G$ an arbitrary {\it
affine} root system and obtain in this way a combinatorial model
for the quantum cohomology ring of the infinite dimensional flag
manifold $LK/T$,  which is investigated in \cite{Ma3}.


\vspace{0.5cm}

{\bf Acknowledgements.} I would like to thank Lisa Jeffrey for a
careful reading of the manuscript and   several helpful comments.
I am also thankful to the referee for several valuable
suggestions.

\newpage

\section{Definition and presentations of $(\H \otimes \bR[\{q_i\}],\star)$}

Our first concern is to show that the combinatorial quantum
product $\star$ described by  equation (\ref{defstar}) is
well-defined.

\noindent {\it Proof of Theorem \ref{defthm}.}  Let us note that
in fact we can define the product $\star$ on $\bR[\{\lambda_i\},\{
q_i\}]$, as follows:
\begin{equation}\label{fstarg}f\star
g:=\psi(\psi^{-1}(f)\psi^{-1}(g))=(\psi^{-1}f)(\{\Lambda_i\},\{q_i\})(g),\end{equation}
$f,g\in\bR[\{\lambda_i\},\{ q_i\}]$. If $g\in I_W\otimes
\bR[\{q_i\}]$, then the last expression in (\ref{fstarg}) is in
$I_W\otimes \bR[\{q_i\}]$ as well (since the latter is  invariant
under any $\Lambda_j$, $1\le j \le l$). We deduce that $I_W\otimes
\bR[\{q_i\}]$ is an ideal of the ring $(\bR[\{\lambda_i\},
\{q_i\}],\star)$. The quotient of the latter ring by the former
ideal is just $(\H\otimes \bR[\{q_i\}],\star)$. It is commutative,
associative and satisfies the grading  condition $\deg (a\star b)
=\deg a+ \deg b$, because the ring
$(\bR[\{\lambda_i\},\{q_i\}],\star)$ is commutative and
associative, and the operator $\Lambda_i$ defined by
(\ref{Lambda})  satisfies $$\deg \Lambda_i(f) = \deg  f +2,$$ for
any  homogeneous polynomial $f\in \bR[\{\lambda_i\},\{q_i\}]$
(provided that $\deg \lambda_i : =2$, $\deg q_i : =4$).

In order to prove the Frobenius property, we only have to check
that \begin{equation}\label{frob}( [\lambda_i] \star [c_v],
[c_w])=( [c_v], [\lambda_i] \star [c_w])\end{equation} for any
$1\le i\le l$, $v,w\in W$. In turn, (\ref{frob}) follows from the
fact that $$[\lambda_i]\star [c_w]=\Lambda_i([c_w])$$ (see
equation (\ref{Lambdastar}) in the introduction), the definition
(\ref{Lambda}) of $\Lambda_i$ and the equation
$$( \Delta_{s_{\alpha}}[c_v],[c_w])=(
[c_v],\Delta_{s_{\alpha}}[c_w]),$$ $v,w \in W$, $\alpha\in
\Phi^+,$ which is a consequence of (\ref{poincare}) and
(\ref{del}).
 \hfill $\square$

We are interested now in obtaining a presentation of the ring
$(\H\otimes \bR[\{q_i\}], \star)$ in terms of generators and
relations. One way\footnote{I am grateful to the referee for
suggesting  me this idea.} of obtaining this is as follows:

\begin{proposition}\label{psiminus} Let $I_W^q$ be the ideal
of $\bR[\{\lambda_i\},\{q_i\}]$ generated by
$\psi^{-1}(u_1),\ldots, \psi^{-1}(u_l)$. The map $$\psi^{-1}
:(\H\otimes \bR[\{q_i\}]=\bR[\{\lambda_i\},\{q_i\}]/(I_W\otimes
\bR[\{q_i\}]),\star )\to \bR[\{\lambda_i\}, \{q_i\}]/I_W^q$$ given
by
$$f \ {\rm mod} \ I_W\otimes \bR[\{q_i\}] \  \mapsto \ \psi^{-1}(f) \ {\rm mod} \
I_W^q,$$ $f\in \bR[\{\lambda_i\},\{q_i\}]$ is a ring isomorphism.
\end{proposition}

\begin{proof} From the definition (\ref{fstarg}) we can see that
\begin{equation}\label{ringisom}\psi^{-1} :
(\bR[\{\lambda_i\},\{q_i\}],\star) \to
(\bR[\{\lambda_i\},\{q_i\}],\cdot)
\end{equation} is a ring
isomorphism. As pointed out before (see the proof of Theorem
\ref{defthm}), the combinatorial quantum cohomology ring
$(\H\otimes \bR[\{q_i\}],\star)$ is the quotient of the ring
$(\bR[\{\lambda_i\},\{q_i\}],\star)$ by its ideal $I_W\otimes
\bR[\{q_i\}]$. Note that the latter --- regarded as an ideal of
$(\bR[\{\lambda_i\},\{q_i\}],\star)$ --- is generated by the same
fundamental $W$-invariant polynomials $u_1,\ldots, u_l$. This is
because for any $f\in \bR[\{\lambda_i\},\{q_i\}]$ we have
$$f\star u_k = \psi^{-1}(f)\cdot u_k,$$
$k=1,\ldots, l$, and $\psi^{-1}$ is a bijective map. Consequently,
the ring isomorphism (\ref{ringisom}) maps the quotient of
$(\bR[\{\lambda_i\},\{q_i\}],\star)$ by the ideal generated by
$u_1,\ldots, u_l$ isomorphically onto the quotient of
$(\bR[\{\lambda_i\},\{q_i\}],\cdot)$ by the ideal generated by
$\psi^{-1}(u_1),\ldots, \psi^{-1}(u_l)$.
\end{proof}

As pointed out out in the introduction, we are also able to deduce
B. Kim's presentation [Kim] for the combinatorial quantum
cohomology ring. In fact
 Theorem \ref{firsttheorem} is a straightforward consequence of the
 following result, which was proved in \cite{Ma2}:

\begin{theorem}\label{gen-rel}{\rm (\cite{Ma2})}
Let $\bullet$ be an $\bR[\{q_i\}]$-bilinear product on $\H\otimes
\bR[\{q_i\}]$ with the following properties:
\begin{itemize}
\item[(i)] $\bullet$ is commutative
\item[(ii)] $\bullet$ is associative
\item[(iii)] $\bullet$ is a
deformation of the usual product, in the sense that if we formally
replace all $q_i$ by $0$, we obtain the usual product on $\H$
\item[(iv)] $(\H\otimes \bR[\{q_i\}], \bullet)$ is a graded ring with respect to $\deg
[\lambda_i] =2$ and $\deg q_i=4$
\item[(v)] $[\lambda_i]\bullet [\lambda_j] =[\lambda_i][\lambda_j]
+\delta_{ij}q_j$
\item[(vi)] $d_i([\lambda_j] \bullet a)_d=d_j([\lambda_i]
\bullet a)_d$, for any $a\in \H$, $1\leq i,j\leq l$, and
$d=(d_1,\ldots, d_l)\ge 0$ (here we use the notation
$[\lambda_i]\bullet a=\sum_{d=(d_1,\ldots, d_l)\ge
0}([\lambda_i]\bullet a)_dq_1^{d_1}\ldots q_l^{d_l},$ with
$([\lambda_i]\bullet a)_d\in \H$).
\end{itemize}
Then the following relation holds in the ring
$(\H\otimes\bR[\{q_i\}],\bullet)$:
\begin{equation}\label{toda}F_k(\{[\lambda_i]\bullet\},\{-\langle
\alpha_i^{\vee},\alpha_i^{\vee}\rangle q_i\} )= 0,
\end{equation} $1\leq k\leq l,$ where $F_k$ are the integrals of
motion of the Toda lattice associated to the coroot system of $G$
(see below).
 Moreover, the ring $(\H \otimes
\bR[\{q_i\}], \bullet)$ is isomorphic to $\bR[\{
\lambda_i\},\{q_i\}]$ modulo the ideal generated by
 $F_k(\{\lambda_i\},\{ -\langle \alpha_i^{\vee},\alpha_i^{\vee}\rangle q_i\} )$,
$1\leq k \leq l$.
\end{theorem}

The Toda lattice we are referring to in the theorem is the
Hamiltonian system whose phase space is
$(\bR^{2l},\sum_{i=1}^ldr_i\wedge ds_i)$ and Hamiltonian function
$$E=\sum_{i,j=1}^l \langle \alpha_i^{\vee},
\alpha_j^{\vee}\rangle r_i r_j +\sum_{i=1}^l e^{2s_i}.$$ It turns
out (see for instance [Go-Wa]) that this system admits $l$
independent integrals of motion $E=F_1,F_2,\ldots, F_l$, which are
all polynomial functions in variables $r_1,\ldots,
r_l,e^{2s_1},\ldots, e^{2s_l}$ and satisfy the condition
\begin{equation}\label{qequalszero}F_k(\lambda_1,\ldots,\lambda_l,0,\ldots,0)=u_k(\lambda_1,\ldots,\lambda_l),
\end{equation}
where $u_1,\ldots,u_l$ are the fundamental $W$-invariant
polynomials (see section 1). According to Theorem \ref{gen-rel},
the ring $(\H \otimes \bR[\{q_i\}], \bullet)$ is generated by
$[\lambda_1],\ldots,[\lambda_l],q_1,\ldots, q_l$, and the
relations are obtained by taking all polynomials $F_k$ and for
each of them making the replacements
$$r_i\mapsto [\lambda_i]\bullet,\quad e^{2s_i}\mapsto
-\langle \alpha_i^{\vee},\alpha_i^{\vee}\rangle q_i,\quad 1\le
i\le l.$$

It is easy to see that the combinatorial quantum product $\star$
satisfies the hypotheses (i)-(iv) of Theorem \ref{gen-rel}. We
prove condition (v) as follows:
$$[\lambda_i]\star [\lambda_j] =
[\psi(\lambda_i)]\star[\psi(\lambda_j)] =
[\psi(\lambda_i\lambda_j)] = [\Lambda_i(\lambda_j)] =
[\lambda_i\lambda_j +\delta_{ij}q_j],$$ $1\le i,j\le l$.
 In order to prove (vi), we note
that the coefficient of $q^{\alpha^{\vee}}$ in $$[\lambda_j]\star
a=\Lambda_j(a)$$ is
$\lambda_i(\alpha^{\vee})\Delta_{s_{\alpha}}(a)$; thus for the
multi-index
$d=\alpha^{\vee}=\lambda_1(\alpha^{\vee})\alpha_1^{\vee}+\ldots +
\lambda_l(\alpha^{\vee})\alpha_l^{\vee}$ we have
$$d_i([\lambda_j]\star
a)_d=\lambda_i(\alpha^{\vee})\lambda_j(\alpha^{\vee})\Delta_{s_{\alpha}}(a),$$
which is symmetric in $i$ and $j$.

Our next goal is to show that the ``quantum BGG-polynomials" (see
Theorem \ref{bernstein}) $\psi^{-1}(c_w)$, $w\in W$, satisfy a
certain orthogonality property, which was already pointed out for
$G=SL(n,\bC)$  by Fomin, Gelfand and Postnikov in [Fo-Ge-Po]. For
any $f\in \bR[\{\lambda_i\},\{q_i\}]$ we denote by $[f]_q$ its
class modulo $I_W^q$. By Theorem \ref{firsttheorem}, the set
$\{[\psi^{-1}({c}_w)]_q ~|~ w\in W\}$ is a basis of
$\bR[\{\lambda_i\},\{q_i\}]/I_W^q$ as an $\bR[\{q_i\}]$-module.
Define
$$\la [f]_q\ra=\alpha_{w_0}$$
where the elements $\alpha_w$ of $\bR[q_i]$ are defined by
$$[f]_q=\sum_{w\in W}\alpha_w[\psi^{-1}({c}_w)]_q.$$
Consider the pairing $\la \ , \ \ra$ on
$\bR[\{\lambda_i\},\{q_i\}]/I_W^q$ given by
$$\la [f]_q ,[g]_q  \ra =\la [fg]_q   \ra.$$
\begin{proposition}\label{prop} We have that
$$\la [\psi^{-1}({c}_u)]_q ,[\psi^{-1}({c}_v)]_q  \ra =
\left\{
\begin{array}{ll}
                         1,\ {\rm if }\ u=w_0v  \\
                         0, \ {\rm otherwise}
                             \end{array}
                      \right.$$
                      \end{proposition}
\begin{proof}
Write $$[\psi^{-1}({c}_u)\psi^{-1}({c}_v)]_q=\sum_{w\in
W}\alpha_w[\psi^{-1}({c}_w)]_q,$$ which means that the polynomial
\begin{equation}\label{ex}\psi^{-1}({c}_u)\psi^{-1}({c}_v)-\sum_{w\in
W}\alpha_w\psi^{-1}({c}_w)\end{equation} is in $I_W^q$. Consider
 $\psi$ of the expression (\ref{ex}), take into account that
$\psi^{-1}({c}_w)(\{[\lambda_i]\star\},\{q_i\})=[c_w]$ and that
$\psi(I_W^q)=I_W\otimes \bR[\{q_i\}]$ (see Proposition
\ref{psiminus}) and obtain in this way the following equality in
$\H\otimes \bR[\{q_i\}]$:
$$[c_u]\star [c_v]=\sum_{w\in W}\alpha_w[c_w]$$
If $( \ , \  ) $ denotes the usual Poincar\'e
pairing\footnote{Actually its $\bR[\{q_i\}]$-linear extension.}
 on
$\H\otimes \bR[\{q_i\}]$, we deduce that
$$\alpha_{w_0}=( [c_u]\star [c_v], 1)=( [c_u],
[c_v])$$ where we have used the Frobenius property of $\star$. The
orthogonality relation stated in the lemma is a direct consequence
of equation (\ref{poincare}).
\end{proof}

\section{Commutativity of the operators $\Lambda_1,\ldots, \Lambda_l$}

The goal of this section is to provide a proof of Lemma
\ref{peterson}. Let us start with the following recursive
construction  of the elements of $\tilde{\Phi}^+$ (the latter has
been defined immediately after Lemma \ref{inequality}).

\begin{lemma}\label{description} A positive root $\alpha$ is in
$\tilde{\Phi}^+$ if and only if it is simple, or else there exist
 $k\ge 2$ and  $i_1,\ldots, i_k \in\{1,\ldots, l\}$ such that
$$\alpha = s_{i_k}\ldots s_{i_2}(\alpha_{i_1})$$ and
$$\alpha_{i_{j+1}} (s_{i_j} \ldots
s_{i_2}(\alpha_{i_1})^{\vee})=-1,$$ for all $1\le j \le k-1$.
Moreover, the expression $$s_{\alpha} = s_{i_k} \ldots
s_{i_2}s_{i_1}s_{i_2} \ldots s_{i_k}$$ is reduced and we have
$$\alpha^{\vee} = \alpha_{i_1}^{\vee} +\ldots +
\alpha_{i_k}^{\vee},$$ hence $\h(\alpha^{\vee})=k$. All roots
$s_{i_j}\ldots s_{i_2}(\alpha_{i_1})$, $1\le j\le k$, are in
$\tilde{\Phi}^+$.
\end{lemma}

\begin{proof} First  we use induction on $k\ge 1$ to prove that
any root of the form described in the lemma is in
$\tilde{\Phi}^+$. Since any simple root is in $\tilde{\Phi}^+$, we
only have to perform the induction step. Assume that $k\ge 2$. The root
$$\beta:=s_{i_{k-1}}\ldots s_{i_2}(\alpha_{i_1})$$
satisfies the hypotheses of the lemma, hence it is in
$\tilde{\Phi}^+$. Moreover, we have
$\alpha_{i_k}(\beta^{\vee})=-1$, hence
$$\alpha^{\vee}=s_{i_k}(\beta^{\vee}) = \beta^{\vee}+
\alpha_{i_k}^{\vee},$$ which implies that
\begin{equation}\label{alphabeta} \h(\alpha^{\vee}) =
\h(\beta^{\vee})
+1.\end{equation} In particular, $\alpha$ is not a simple root.
Also because $\alpha_{i_k}(\alpha^{\vee}) = 1$, we deduce that the
roots
$$s_{\alpha}(\alpha_{i_k}) = \alpha_{i_k}-\alpha_{i_k}(\alpha^{\vee})\alpha \ {\rm and } \ s_{i_k}s_{\alpha}(\alpha_{i_k})=(\alpha(\alpha_{i_k}^{\vee})\alpha_{i_k}(\alpha^{\vee})-
1)\alpha_{i_k}- \alpha_{i_k}(\alpha^{\vee})\alpha$$   are both
negative. Consequently we have
$$l(s_{\alpha}) = l(s_{i_k}s_{\alpha}s_{i_k}) +2 = l(s_{\beta}) +2
= 2\h(\beta^{\vee}) -1 +2 =2 \h(\alpha^{\vee}) -1,$$ where we have
used (\ref{alphabeta}). Hence $\alpha \in \tilde{\Phi}^+$.

Now we will use induction on $l(s_{\alpha})$ in order to prove
that any element of $\tilde{\Phi}^+$ can be realized in this way.
If $l(s_{\alpha})=1$, then $\alpha$ is simple, hence it is of the
type indicated in the lemma. Assume now that $\alpha\in
\tilde{\Phi}^+$ is not simple. There exists a simple root
$\alpha_i$ such that $\alpha(\alpha_i^{\vee})>0$ (otherwise we
would be led to $\alpha(\alpha^{\vee})\leq 0$). Also
$\alpha_i(\alpha^{\vee})$ must be  strictly positive, hence the
roots
$$s_{\alpha}(\alpha_{i}) = \alpha_{i}-\alpha_{i}(\alpha^{\vee})\alpha \ {\rm and } \ s_{i}s_{\alpha}(\alpha_i)=(\alpha(\alpha_i^{\vee})\alpha_i(\alpha^{\vee})-1)\alpha_i-
\alpha_i(\alpha^{\vee})\alpha$$ are both  negative. We deduce that
$l(s_{i}s_{\alpha}s_{i})=l(s_{\alpha})-2$. From
$$s_{i}
(\alpha)^{\vee}=s_{i}(\alpha^{\vee})=\alpha^{\vee}-\alpha_i(\alpha^{\vee})\alpha_i^{\vee}$$
it follows that $s_i(\alpha)$ is a positive root which satisfies
$\h(s_{i}
(\alpha)^{\vee})=\h(\alpha^{\vee})-\alpha_i(\alpha^{\vee})$. By
Lemma \ref{inequality}, we have that:
$$l(s_{\alpha})=l(s_{i}s_{\alpha}s_{i})+2\leq
2\h(s_{i}(\alpha)^{\vee})-1+2=2\h(\alpha^{\vee})-1+2(1-\alpha_i(\alpha^{\vee}))\leq
2\h(\alpha^{\vee})-1.$$ Since $\alpha\in \tilde{\Phi}^+$, the two
inequalities from the last equation must be equalities. In other
words, $s_i\alpha \in \tilde{\Phi}^+$ and $\alpha_i(\alpha^{\vee})
= 1$, the latter being equivalent to
$\alpha_i((s_i\alpha)^{\vee})=-1$. We  use the induction
hypothesis for $s_i\alpha$, which has the property that
$l(s_{s_i\alpha}) = l(s_is_{\alpha}s_i) =l(s_{\alpha}) -2$ and the
induction step is accomplished.
\end{proof}

The following  property of $\tilde{\Phi}^+$ will be needed later.

\begin{lemma}\label{firstlemma1} If $\alpha, \beta\in
\tilde{\Phi}^+$ are such that
$$l(s_{\alpha}s_{\beta})=l(s_{\alpha})+l(s_{\beta})$$ and
$s_{\alpha}s_{\beta}\neq s_{\beta}s_{\alpha}$, then
$\alpha(\beta^{\vee})<0$.
\end{lemma}

\begin{proof} We use induction on $l(s_{\beta})$. If $\beta$ is
simple, the condition $l(s_{\alpha}s_{\beta})=l(s_{\alpha})+1$ is
equivalent to the fact that the root
$s_{\alpha}(\beta)=\beta-\beta(\alpha^{\vee})\alpha$ is positive,
which implies $\beta(\alpha^{\vee})\leq 0$, and then
$\alpha(\beta^{\vee})\leq 0$. We cannot have $\alpha(\beta^{\vee})=0$,
since otherwise $s_{\alpha}$ and $s_{\beta}$ would commute.

The induction step will follow now. Let us assume first that the
root system involved here is not of type $G_2$. Consider  $\alpha,
\beta\in \tilde{\Phi}^+$  both non-simple; by Lemma
\ref{description}, $\beta$ is of the form
$\beta=s_i(\tilde{\beta})$, where $\tilde{\beta}\in
\tilde{\Phi}^+$ and $\alpha_i(\tilde{\beta}^{\vee})=-1$. Suppose
that $\alpha(\beta^{\vee})\geq 1$. Since
$\alpha_i(\beta^{\vee})=1$, the root $s_{\beta}(\alpha_i)
=\alpha_i-\beta$ is negative, hence
$l(s_{\beta}s_i)=l(s_{\beta})-1$. From
$l(s_{\alpha}s_{\beta})=l(s_{\alpha})+l(s_{\beta})$, we deduce now
that $l(s_is_{\alpha})=l(s_{\alpha})+1$, hence the root
$s_{\alpha}(\alpha_i)=\alpha_i-\alpha_i(\alpha^{\vee})\alpha$ is
positive, which implies $\alpha_i(\alpha^{\vee})\le 0$.

{\it Claim.}  $\alpha_i(\alpha^{\vee})\ne 0$.

Because otherwise
$s_i$ and $s_{\alpha}$ commute, hence
\begin{align}l(s_{\tilde{\beta}}s_{\alpha})&=l(s_{\tilde{\beta}}s_is_{\alpha}s_i)\nonumber \\{}& =
l(s_{\tilde{\beta}}s_is_{\alpha}) -1 \nonumber
\\{}&=l(s_is_{\tilde{\beta}}s_is_{\alpha}) -2
\nonumber \\{}& =l(s_{{\beta}}s_{\alpha}) -2
\nonumber \\{}&=l(s_{{\beta}}) -2 +l(s_{\alpha})
\nonumber \\{}&=l(s_{\tilde{\beta}})  +l(s_{\alpha})\nonumber
 \end{align}
where the second equality holds since $l(s_{\tilde{\beta}}s_is_{\alpha})=l(s_{\tilde{\beta}})+
l(s_{\alpha})+1 > l(s_{\tilde{\beta}}s_{\alpha})$. By the induction hypothesis, we must have
$\tilde{\beta}(\alpha^{\vee})\le 0$. On the other hand we have
$$\tilde{\beta}(\alpha^{\vee}) = s_i\beta(\alpha^{\vee}) = \beta(s_i\alpha^{\vee})=
\beta(\alpha^{\vee})$$
the last number being strictly positive. This contradiction concludes the claim.

From the claim we deduce  that
\begin{equation}\label{ge}\alpha(\tilde{\beta}^{\vee}) = \alpha(\beta^{\vee}) - \alpha(\alpha_i^{\vee}) \ge 2.\end{equation}
Since the root system is not of type $G_2$, we must have equality in  (\ref{ge}), hence
\begin{equation}\label{minusone}\alpha(\alpha_i^{\vee}) =-1.\end{equation} We distinguish the following two possibilities:

(i) $\alpha \ne \tilde{\beta}.$
From (\ref{ge}) we deduce that $||\tilde{\beta}|| < ||\alpha||$.
Since $||\tilde{\beta}||=||s_i \tilde{\beta}|| = ||\beta||$, we have that
$||\beta|| <||\alpha||$, hence $ \alpha({\beta}^{\vee}) \ge 2$. Consequently,
\begin{equation}\label{ge3}\alpha(\tilde{\beta}^{\vee}) = \alpha(\beta^{\vee}) - \alpha(\alpha_i^{\vee}) \ge 3,\end{equation}
which cannot happen as long as  the root system is not of type $G_2$.

(ii) $\alpha = \tilde{\beta}$. This means that $\beta = s_i(\alpha)$,
\begin{equation}\label{first}\alpha_i(\alpha^{\vee})=-1,\end{equation}
and $\beta^{\vee} = \alpha^{\vee}+\alpha_i^{\vee}$. From (\ref{minusone}) and (\ref{first})
we deduce  that $$s_{\alpha}s_is_{\alpha}(\alpha_i)=-\alpha,$$
which is a negative root, hence
$$l(s_{\alpha}s_{\beta})=l(s_{\alpha}s_is_{\alpha}s_i)=l(s_{\alpha}s_is_{\alpha})-1\le
l(s_{\alpha}) + l(s_is_{\alpha})-1 = l(s_{\alpha}) + l(s_is_{\alpha}s_i)-2=l(s_{\alpha}) +
l(s_{\beta})-2.$$
This is a contradiction.

Now let us consider the case when the
root system is of type $G_2$. Let
$\alpha_1$, $\alpha_2$ be the standard basis of the root system
$G_2$, with $||\alpha_1||
> ||\alpha_2||$. By Lemma \ref{description} we can see
that $\tilde{\Phi}^+$ consists of $\alpha_1$, $\alpha_2$, $s_2(\alpha_1)=\alpha_1+3\alpha_2$,
and $s_1s_2(\alpha_1)=2\alpha_1+3\alpha_2$.
Since $\alpha(\beta^{\vee})\ge 1$ and none of $\alpha$ and $\beta$ is simple, we can only have
 $\alpha = s_1 s_2
(\alpha_1)$ and  $\beta = s_2(\alpha_1)$, which implies
$s_{\alpha}=s_1s_2s_1s_2s_1$ and $s_{\beta}= s_2s_1s_2$,  hence
$s_{\alpha}s_{\beta}=s_1s_2s_1s_2s_1s_2s_1s_2=(s_1s_2)^4$; but the
latter is the same as $(s_1s_2)^2$, having length 4, which is
strictly less than $l(s_{\alpha}) + l(s_{\beta}) = 5+3=8$. The
contradiction shows that also in this case we must have
$\alpha(\beta^{\vee})<0$.

\end{proof}

\begin{lemma}\label{elementof} If $\alpha, \beta \in \tilde{\Phi}^+$ with
$l(s_{\alpha}s_{\beta}) = l(s_{\alpha})+l(s_{\beta})$ and
$s_{\alpha}s_{\beta}\ne s_{\beta}s_{\alpha}$, then there exists
$\gamma\in \tilde{\Phi}^+$ such that $$ \alpha^{\vee}+\beta^{\vee}
=
\gamma^{\vee}.$$ 
\end{lemma}

\begin{proof} By Lemma \ref{firstlemma1}, one of the numbers $\alpha(\beta^{\vee})$
and $\beta(\alpha^{\vee})$ is $-1$.   We will actually prove that
if $\beta(\alpha^{\vee})=-1$ then $s_{\beta}(\alpha)\in
\tilde{\Phi}^+$ (it is obvious that
$s_{\beta}(\alpha)^{\vee}=s_{\beta}(\alpha^{\vee}) = \alpha^{\vee}
+\beta^{\vee}$). We will use induction on $l(s_{\beta})$. If
$\beta$ is simple, the result follows immediately from Lemma
\ref{description}. Consider now the case when $\beta\in
\tilde{\Phi}^+$ is non-simple; by Lemma \ref{description}, $\beta$
is of the form $\beta=s_i(\tilde{\beta})$, where $\tilde{\beta}\in
\tilde{\Phi}^+$ and $\alpha_i(\tilde{\beta}^{\vee})=-1$. From
$l(s_{\beta}s_i)=l(s_{\beta})-1$ and $l(s_{\alpha}s_{\beta})=
l(s_{\alpha})+l(s_{\beta})$ it follows that $l(s_{\alpha}s_i)=
l(s_{\alpha})+1$, hence $s_{\alpha}(\alpha_i)=\alpha_i
-\alpha_i(\alpha^{\vee})\alpha$ is positive, which means
$\alpha_i(\alpha^{\vee})\leq 0$. We show that the only possible
values for $\alpha_i(\alpha^{\vee})$ are  $-1$ and $0$. Otherwise,
the root system is {\it not} simply laced and the roots $\alpha$
and $\alpha_i$ are
 short, respectively long; on the other hand,
$\alpha_i(\beta^{\vee})=1$, so $||\alpha_i||\leq ||\beta||$
and $\beta(\alpha^{\vee})=-1$, so $||\beta|| \leq ||\alpha||$,
which gives a contradiction.

{\it Case 1.} $\alpha_i(\alpha^{\vee})=0$. This implies
$s_i(\alpha)=\alpha$, hence
$$-1=\beta(\alpha^{\vee})=s_i(\tilde{\beta})(\alpha^{\vee})=
\tilde{\beta}(s_i(\alpha)^{\vee})=\tilde{\beta}(\alpha^{\vee}).$$
From the induction hypothesis,
$s_{\tilde{\beta}}(\alpha)=s_is_{\beta}(\alpha):=\gamma$ is in
$\tilde{\Phi}^+$. We also have that
$$\alpha_i(\gamma^{\vee})=\alpha_i (s_is_{\beta}(\alpha^{\vee}))=
-\alpha_i(s_{\beta}(\alpha^{\vee}))=-\alpha_i(\alpha^{\vee}
+\beta^{\vee})=-1.$$ By Lemma \ref{description}, the root
$s_i(\gamma)=s_{\beta}(\alpha)$ is in $\tilde{\Phi}^+$.

{\it Case 2.} $\alpha_i(\alpha^{\vee})=-1$. We have again that
$$-1=\beta(\alpha^{\vee})=s_i(\tilde{\beta})(\alpha^{\vee})=
\tilde{\beta}(s_i(\alpha)^{\vee}).$$ By Lemma \ref{description},
the root $s_i(\alpha)$ is in $\tilde{\Phi}^+$, and from the
induction hypothesis we deduce that
$s_{\tilde{\beta}}(s_i(\alpha))=s_is_{\beta}(\alpha):=\gamma$ is
also in $\tilde{\Phi}^+$. But, as before,
$$\alpha_i(\gamma^{\vee})=-\alpha_i(\alpha^{\vee}+\beta^{\vee}),$$
the right hand side being now $0$. We deduce that
$\gamma=s_i(\gamma)= s_{\beta}(\alpha)$.\end{proof}

We are now able to prove Lemma \ref{peterson}:

{\it Proof of Lemma  \ref{peterson}}. Denote by $\lambda_i^*$ the
operator of  multiplication by $\lambda_i$ on
$\bR[\{\lambda_1,\ldots, \lambda_l\}]$, $1\le i\le l$. The
following formula can be found for instance in [Hi, Ch. IV,
section 3]:
\begin{equation}\label{hiller}\Delta_w\lambda_i^*-w\lambda_i^*
w^{-1}\Delta_w= \sum_{\beta\in
\Phi^+,l(ws_{\beta})=l(w)-1}\lambda_i(\beta^{\vee})
\Delta_{ws_\beta},\end{equation} where $w\in W$. Put
$w=s_{\alpha}$ in (\ref{hiller}) and obtain that: $$
\Delta_{s_{\alpha}}\lambda_i^*=(\lambda_i^*
-\lambda_i(\alpha^{\vee})\alpha^*)\Delta_{s_{\alpha}}
+\sum_{\gamma\in \Phi^+,
l(s_{\alpha}s_{\gamma})=l(s_{\alpha})-1}\lambda_i(\gamma^{\vee})
\Delta_{s_{\alpha}s_{\gamma}}.$$ We deduce that:
\begin{align} \Lambda_j\Lambda_i=&(\lambda_j\lambda_i)^*
+\sum_{\alpha\in\tilde{\Phi}^+}\lambda_i(\alpha^{\vee})
q^{\alpha^{\vee}}\lambda_j^*\Delta_{s_{\alpha}}+
\sum_{\alpha\in\tilde{\Phi}^+}\lambda_j(\alpha^{\vee})
q^{\alpha^{\vee}}\lambda_i^*\Delta_{s_{\alpha}}\nonumber\\{}&
-\sum_{\alpha\in\tilde{\Phi}^+}\lambda_j(\alpha^{\vee})
\lambda_i(\alpha^{\vee})q^{\alpha^{\vee}}\alpha^*\Delta_{s_{\alpha}}\nonumber\\{}&
+\sum_{\alpha\in\tilde{\Phi}^+,\gamma\in \Phi^+,
l(s_{\alpha}s_{\gamma})= l(s_{\alpha})-1}\lambda_j(\alpha^{\vee})
\lambda_i(\gamma^{\vee})q^{\alpha^{\vee}}\Delta_{s_{\alpha}
s_{\gamma}}\nonumber \\{}& + \sum_{\beta,\delta \in
\tilde{\Phi}^+, l(s_{\beta}s_{\delta})=l(s_{\beta})+l(s_{\delta})}
\lambda_j(\beta^{\vee})
\lambda_i(\delta^{\vee})q^{\beta^{\vee}+\delta^{\vee}}
\Delta_{s_{\beta}s_{\delta}}\nonumber\\{} =&(\lambda_j\lambda_i)^*
+\sum_{\alpha\in\tilde{\Phi}^+}\lambda_i(\alpha^{\vee})
q^{\alpha^{\vee}}\lambda_j^*\Delta_{s_{\alpha}}+\sum_{\alpha\in\tilde{\Phi}^+}\lambda_j(\alpha^{\vee})
q^{\alpha^{\vee}}\lambda_i^*\Delta_{s_{\alpha}}\nonumber\\
{}& -\sum_{\alpha\in\tilde{\Phi}^+}\lambda_j(\alpha^{\vee})
\lambda_i(\alpha^{\vee})q^{\alpha^{\vee}}\alpha^*\Delta_{s_{\alpha}}\nonumber\\
{}& +\sum_{\alpha \ {\rm
simple}}\lambda_j(\alpha^{\vee})\lambda_i(\alpha^{\vee})q^{\alpha^{\vee}}\nonumber\\
{}& +\sum_{\beta,\delta \in \tilde{\Phi}^+,
l(s_{\beta}s_{\delta})=l(s_{\beta})+l(s_{\delta}),
s_{\beta}s_{\delta}=s_{\delta}s_{\beta}}\lambda_j(\beta^{\vee})
\lambda_i(\delta^{\vee})q^{\beta^{\vee}+\delta^{\vee}}
\Delta_{s_{\beta}s_{\delta}}\nonumber\\
{}&+\sum_{\alpha\in\tilde{\Phi}^+, \gamma\in \Phi^+,
l(s_{\alpha}s_{\gamma})= l(s_{\alpha})-1\ge
1}\lambda_j(\alpha^{\vee})
\lambda_i(\gamma^{\vee})q^{\alpha^{\vee}}\Delta_{s_{\alpha}
s_{\gamma}}\nonumber \\{}& + \sum_{\beta,\delta \in
\tilde{\Phi}^+,
l(s_{\beta}s_{\delta})=l(s_{\beta})+l(s_{\delta}),s_{\beta}s_{\delta}\ne
s_{\delta}s_{\beta}} \lambda_j(\beta^{\vee})
\lambda_i(\delta^{\vee})q^{\beta^{\vee}+\delta^{\vee}}
\Delta_{s_{\beta}s_{\delta}}\nonumber
\end{align}

Denote by $\Sigma_{ij}$ the sum of the last two sums: the rest is
obviously invariant under the operation of interchanging $i$ and $
j$.

We will show that $\Sigma_{ij}$ is symmetric in $i$ and $j$. To
this end, let us take first two arbitrary elements $\beta, \delta$
of $\tilde{\Phi}^+$ with
$l(s_{\beta}s_{\delta})=l(s_{\beta})+l(s_{\delta})$ and
$s_{\beta}s_{\delta}\neq s_{\delta}s_{\beta}$; by Lemma
\ref{firstlemma1} and Lemma \ref{elementof}, there exists
$\alpha\in \tilde{\Phi}^+$ such that
 $\alpha^{\vee}=\beta^{\vee}+\delta^{\vee}$; we will show
that: \begin{itemize} \item there exists a unique $\gamma \in
\Phi^+$ with $s_{\alpha}s_{\gamma}=s_{\beta}s_{\delta} \ {\rm and
} \ l(s_{\alpha}s_{\gamma})=l(s_{\alpha})-1,$ \item for $\gamma$
determined above, the sum
$$ \lambda_j(\alpha^{\vee})\lambda_i(\gamma^{\vee})
\Delta_{s_{\alpha}s_{\gamma}} +
\lambda_j(\beta^{\vee})\lambda_i(\delta^{\vee})\Delta_{s_{\beta}s_{\delta}}
=(\lambda_j(\alpha^{\vee})\lambda_i(\gamma^{\vee}) +
\lambda_j(\beta^{\vee})\lambda_i(\delta^{\vee}))
\Delta_{s_{\beta}s_{\delta}} := S_{ij}^{\beta,\delta}
\Delta_{s_{\beta}s_{\delta}}$$ is symmetric in $i$ and $j$.
\end{itemize}
 By Lemma \ref{firstlemma1}, we
distinguish the following two cases:

{\it Case 1.} $\beta(\delta^{\vee})=-1$, which implies
$\alpha=s_{\beta}(\delta)$, so the condition
$s_{\alpha}s_{\gamma}=s_{\beta}s_{\delta}$ is equivalent
$\gamma=\beta$. Note that
$$l(s_{\alpha})=2\h(\alpha^{\vee})-1=2(\h(\beta^{\vee})+ \h(\delta^{\vee}))-1=
l(s_{\beta}s_{\delta})+1=l(s_{\alpha}s_{\gamma})+1.$$ We deduce
that
$$S_{ij}^{\beta,\delta} =\lambda_j(\alpha^{\vee})\lambda_i(\beta^{\vee})+
\lambda_j(\beta^{\vee})\lambda_i(\delta^{\vee})=
\lambda_j(\beta^{\vee})\lambda_i(\beta^{\vee})
+\lambda_j(\delta^{\vee})\lambda_i(\beta^{\vee})
+\lambda_j(\beta^{\vee})\lambda_i(\delta^{\vee})$$ which is
obviously symmetric in $i$ and $j$.

{\it Case 2.} $\delta(\beta^{\vee})=-1$, which implies that
$\alpha=s_{\delta}(\beta)$, so this time the condition
$s_{\alpha}s_{\gamma}=s_{\beta}s_{\delta}$ is equivalent to
$\gamma=\pm s_{\alpha}(\delta)$. Because
$\delta(\alpha^{\vee})=1$, the number $\alpha(\delta^{\vee})$ is
strictly positive, hence the root $s_{\alpha}(\delta)^{\vee} =
s_{\alpha}(\delta^{\vee})=
\delta^{\vee}-\alpha(\delta^{\vee})\alpha^{\vee} =
\delta^{\vee}-\alpha(\delta^{\vee})(\beta^{\vee} + \delta^{\vee})$
is negative, so we must have $\gamma =-s_{\alpha}(\delta)$. We
have again that
$$l(s_{\alpha})=2\h(\alpha^{\vee})-1=2(\h(\beta^{\vee})+
\h(\delta^{\vee}))-1=
l(s_{\beta}s_{\delta})+1=l(s_{\alpha}s_{\gamma})+1.$$
 This time we can express $S_{ij}^{\beta,\delta}$ as follows:
\begin{align}S_{ij}^{\beta,\delta} &=-\lambda_j(\alpha^{\vee})\lambda_i(s_{\alpha}(\delta^{\vee}))+
\lambda_j(\beta^{\vee})\lambda_i(\delta^{\vee})
 \nonumber\\{}&=-\lambda_j(\alpha^{\vee})
(\lambda_i(\delta^{\vee})-\lambda_i(\alpha^{\vee})\alpha(\delta^{\vee}))
+\lambda_j(\beta^{\vee})\lambda_i(\delta^{\vee})
\nonumber\\{}&=-\lambda_j(\delta^{\vee})\lambda_i(\delta^{\vee})
+\lambda_j(\alpha^{\vee})\lambda_i(\alpha^{\vee})\alpha(\delta^{\vee}),\nonumber\end{align}
 which is again symmetric in
$i$ and $j$.

In order to complete the proof, we must
 take $\alpha\in \tilde{\Phi}^+$ non-simple, $\gamma \in \Phi^+$
with $l(s_{\alpha}s_{\gamma})=l(s_{\alpha})-1$ and show that there
exists $\beta,\delta \in \tilde{\Phi}^+$ with
$\beta^{\vee}+\delta^{\vee}= \alpha^{\vee}$,
$s_{\alpha}s_{\gamma}=s_{\beta}s_{\delta}$ and
$s_{\beta}s_{\delta}\neq s_{\delta}s_{\beta}$. Consider the
reduced decomposition $s_{\alpha}=s_{i_k} \ldots
s_{i_2}s_{i_1}s_{i_2} \ldots s_{i_k}$ given by Lemma
\ref{description}. By the ``strong exchange condition" (see for
instance [Hu, section 5.8]) we distinguish the following two
cases:

{\it Case A.} $s_{\gamma}= s_{i_k} \ldots
s_{i_{j+1}}s_{i_j}s_{i_{j+1}} \ldots s_{i_k}$ for some $j$ between
$2$ and $k$. We deduce that $\gamma=s_{i_k} \ldots
s_{i_{j+1}}(\alpha_{i_j})$, the latter being a positive root since
the expression $s_{i_k} \ldots s_{i_{j+1}}s_{i_j}$ is reduced. We
notice that
$$\gamma(\alpha^{\vee})=s_{i_k} \ldots s_{i_{j+1}}(\alpha_{i_j})
(s_{i_k} \ldots s_{i_{2}}(\alpha_{i_1}^{\vee}))=
\alpha_{i_j}(s_{i_j} \ldots s_{i_{2}}(\alpha_{i_1}^{\vee}))=1,$$
where we have used Lemma \ref{description}.  Set $\beta=\gamma$
and $\delta=s_{\gamma}(\alpha)$, so that $\delta^{\vee}=
s_{\gamma}(\alpha^{\vee})=\alpha^{\vee}-\gamma^{\vee}$, which
implies $\alpha^{\vee}=\beta^{\vee}+\delta^{\vee}$. We obviously
have $s_{\beta}s_{\delta}=s_{\alpha}s_{\gamma}$, hence
$$l(s_{\beta}s_{\delta})=2\h(\alpha^{\vee})-2= 2h(\beta^{\vee}) -1 +
2\h(\gamma^{\vee})-1.$$ From Lemma \ref{inequality} we deduce that
$\beta$ and $\delta$ are both in $\tilde{\Phi}^+$ and
$l(s_{\beta}s_{\delta}) =l(s_{\beta})+l(s_{\delta})$.

{\it Case B.} \begin{align} s_{\gamma}&=s_{i_k} \ldots
s_{i_2}s_{i_1}s_{i_2} \ldots s_{i_{j-1}}s_{i_j}s_{i_{j-1}} \ldots
s_{i_2}s_{i_1}s_{i_2} \ldots s_{i_k}\nonumber
\\{}&=s_{\alpha}s_{i_k}\ldots s_{i_{j+1}}s_{i_j}s_{i_{j+1}} \ldots
s_{i_k} s_{\alpha}\nonumber\end{align} which implies that
\begin{eqnarray*} \gamma=-s_{\alpha}(s_{i_k} \ldots
s_{i_{j+1}}(\alpha_{i_j}))= s_{i_k} \ldots s_{i_2}s_{i_1}s_{i_2}
\ldots s_{i_{j-1}}(\alpha_{i_j}).\end{eqnarray*} A straightforward
calculation shows that   $\gamma(\alpha^{\vee})=1$. We set
$\delta=-s_{\alpha}(\gamma)$, and $\beta=-s_{\alpha}s_{\gamma}
(\alpha)$ (it is not difficult to see that both
$s_{\alpha}(\gamma)$ and $s_{\alpha}s_{\gamma} (\alpha)$ are
negative roots). We have that $\delta^{\vee}=
-\gamma^{\vee}+\alpha(\gamma^{\vee})\alpha^{\vee}$ and
$\beta^{\vee}=\gamma^{\vee}-(\alpha(\gamma^{\vee})-1)\alpha^{\vee}$,
which implies that $\beta^{\vee}+\delta^{\vee}=\alpha^{\vee}$. We
can easily check that $s_{\beta}s_{\delta}=s_{\alpha} s_{\gamma}$.
As in the previous situation, we show that $\beta$ and $\delta$
are both in in $\tilde{\Phi}^+$ and we have
$l(s_{\beta}s_{\delta})=l(s_{\beta})+l(s_{\delta})$. \hfill
$\square$

\bibliographystyle{abbrv}

\end{document}